\documentclass{article}
\usepackage{amsmath,amssymb,latexsym,amsfonts,enumerate,theorem,longtable,amscd}
\newtheorem{definition}{Definition}[section]\newtheorem{proposition}[definition]{Proposition}
\newtheorem{lemma}[definition]{Lemma}
\begin{document}

\begin{center}
{\bf\large A new $Z_3$-graded quantum group}
\end{center}

\begin{center}
Salih Celik \footnote{E-mail: sacelik@yildiz.edu.tr}

Department of Mathematics, Yildiz Technical University, DAVUTPASA-Esenler, 
Istanbul, 34210 TURKEY. 
\end{center}

\noindent{\bf MSC:} 17B37; 81R60

\noindent{\bf Keywords:} $Z_3$-graded exterior algebra, $Z_3$-graded quantum group, $Z_3$-graded Hopf algebra.

\begin{abstract}
We introduce a $Z_3$-graded version of exterior (Grassmann) algebra with two generators and using this object we obtain a new $Z_3$-graded quantum group denoted by $O(\widetilde{GL}_q(2))$. We also discuss some properties of ${ O}(\widetilde{GL}_q(2))$.
\end{abstract}

\section{Introduction}\label{sec1}

Quantum plane \cite{Manin1} is a well known example in quantum group theory. One specific approach to represent a quantum group is to introduce quantum plane (and its dual). When there exists an appropriate set of noncommuting variables spanning linearly a representation space, the endomorphisms on that space preserving the noncommutative structure allows to set up a quantum group. The natural extension to $Z_2$-graded space was introduced in \cite{Manin2}.
The present work starts a $Z_3$-graded version of the exterior plane, denoted by $\tilde{\mathbb R}_q^{0|2}$, where $q$ is a cubic root of unity. In this case, of course, it will not go back to the original objects. The term "plane" is used as a formal title based upon its construction. Following the approach of the Manin's to quantum group $GL(2)$ we see that there exists a $Z_3$-graded (quantum) group acting on the $Z_3$-graded exterior plane. A detailed discussion of this group are given in Sect. 3. In \cite{Chung} Chung finds commutation relations between the elements of a $Z_3$-graded quantum 2x2 matrix using the differential schema established on quantum (1+1)-superplane. With a similar idea, in \cite{Celik} the author obtains similar (but not all the same) relations. However, all structures introduced in the present study are completely different from both \cite{Chung} and \cite{Celik} except for matrix $T$.

\section{$Z_3$-graded planes}\label{sec2}

The aim of this section is to introduce the $Z_3$-graded version of the exterior algebra and its dual. It is known that the Manin's quantum plane is introduced as a $q$-deformation of commutative plane in the sense that it becomes the classical plane when $q$ is equal to 1. In our case, the parameter $q$ is a cubic root of unity and there is no return. To understand what this means, let's begin with recalling some facts about the exterior algebra.

\subsection{$Z_3$-gradation}\label{sec2.1}

A $Z_3$-graded vector space is a vector space $V$ together with a decomposition $V = V_{\bar{0}}\oplus V_{\bar{1}}\oplus V_{\bar{2}}$. Members of $V_{\bar{0}}\oplus V_{\bar{1}}\oplus V_{\bar{2}}$ are called homogeneous elements. The {\it grade} (or {\it degree}) of a homogenous element $v\in V_i$ is denoted by $\tau(v)=i$, $i\in Z_3$. An element in $V_{\bar{0}}$ (resp. $V_{\bar{1}}$ and $V_{\bar{2}}$) is of degree 0 (resp. 1 and 2).

A $Z_3$-graded algebra ${ A}$ is a $Z_3$-graded vector space ${ A} = { A}_{\bar{0}}\oplus { A}_{\bar{1}}\oplus { A}_{\bar{2}}$ which is also an associative algebra such that ${ A}_i \cdot { A}_j \subset { A}_{i+j}$ or, equivalently, $\tau(\xi_1\cdot\xi_2)=\tau(\xi_1)+\tau(\xi_2)$ for all homogeneous elements $\xi_1,\xi_2 \in { A}$.

\subsection{The algebra of functions on the $Z_3$-graded exterior plane}\label{sec2.2}

A possible way to generalize the $Z_3$-graded exterior plane is to increase the power of nilpotency of its generators and to impose a $Z_3$-graded
commutation relation on the generators. We will assume that $q$ is a cubic root of unity.

It is needed to put the wedge product between the coordinates of exterior plane, but it does not matter in the $Z_3$-graded case.

\begin{definition} 
Let $O(\tilde{\mathbb R}_q^{0|2})$ be the algebra with the generators $\theta$ and $\varphi$ obeying the relations
\begin{equation} \label{1}
\theta\cdot\varphi = q^2 \varphi\cdot\theta, \quad \theta^3 = 0 = \varphi^3
\end{equation}
where the coordinates $\theta$ and $\varphi$ are of grade 1 and 2, respectively. We call $O(\tilde{\mathbb R}_q^{0|2})$ the algebra of functions on the Z$_3$-graded exterior plane $\tilde{\mathbb R}_q^{0|2}$.
\end{definition}

\begin{definition}
The $Z_3$-graded plane $\tilde{\mathbb R}_q^{*0|2}$ with the function algebra
$${ O}(\tilde{\mathbb R}_q^{*0|2}) = K\{\xi,x\}/\langle\xi x - x \xi\rangle$$
is called $Z_3$-graded dual exterior plane where the generators $\xi$, $x$ are of degree 2, 0, respectively.
\end{definition}

Hence, in accordance with Definition 2.2, we have
\begin{equation} \label{2}
\tilde{\mathbb R}_q^{*0|2} \ni \begin{pmatrix} \xi \\ x \end{pmatrix} \,\, \Longleftrightarrow \,\, \xi x = x \xi.
\end{equation}

\section{The $Z_3$-graded (quantum) group}\label{sec3}

The algebraic group $SL(2,{\mathbb C})$ has coordinate algebra $O(SL(2,{\mathbb C}))$. This algebra is the quotient of the commutative polynomial
algebra ${\mathbb C}[a,b,c,d]$ by the two-sided ideal generated by the element $ad-bc-1$ where the indeterminates $a,b,c,d$ are the coordinate functions on $SL(2,{\mathbb C})$. Using the group structure in $SL(2,{\mathbb C})$, we can encode it in terms of maps $m$ (multiplication),  $\eta$ (identity) and $S$ (inversion). Dualizing these maps to $O(SL(2,{\mathbb C}))$, we get the corresponding co-maps called comultiplication $\Delta$, counit $\epsilon$, and antipode $S$, respectively. The axioms for the group structure of $SL(2,{\mathbb C})$, in terms of the maps, are then reversed giving us relations among the co-maps. The natural axioms satisfied in $O(SL(2,{\mathbb C}))$ by the maps $m$, $\eta$, $\Delta$, $\epsilon$ and $S$, it makes a Hopf algebra.
The quantum group $O_q(SL(2,{\mathbb C}))$ is a noncommutative deformation of $O(SL(2,{\mathbb C}))$. General concepts related to quantum groups (Hopf algebras) can be found in the books of Klimyk and Schm\"udgen \cite{Kli-sch} or Majid \cite{Majid1}.

In this section, we will consider the 2x2 matrices acting on the $Z_3$-graded exterior plane and will discuss the properties of such matrices. So, let
$a$, $\beta$, $\gamma$, $d$ be elements of an algebra ${ A}$ where the generators $a$ and $d$ are of degree 0, the generators $\gamma$ and $\beta$ are of degree 1 and 2, respectively. Let $\tilde{M}(2)$ be defined as the polynomial algebra $k[a,\beta,\gamma,d]$.
It will sometimes be convenient and more illustrative to write a point $(a,\beta,\gamma,d)$ of $\tilde{M}(2)$ in the matrix form
\begin{eqnarray} \label{3}
T = \begin{pmatrix} a & \beta \\ \gamma & d \end{pmatrix} =(t_{ij}).
\end{eqnarray}
We constitute the $Z_3$-graded matrix algebra $\tilde{M}(2)$ as follows: We divided the algebra $\tilde{M}(2)$ into three parts in form
$\tilde{M}(2)={ A}_{\bar{0}}\oplus { A}_{\bar{1}}\oplus { A}_{\bar{2}}$. In this case, if a matrix has the form of
\begin{eqnarray*}
T_0 = \begin{pmatrix} a & 0 \\ 0 & d \end{pmatrix}, \quad (\mbox{resp.} \quad T_1 = \begin{pmatrix} 0 & 0 \\ \gamma & 0 \end{pmatrix}, \quad T_2 = \begin{pmatrix} 0 & \beta \\ 0 & 0 \end{pmatrix}),
\end{eqnarray*}
then it is an element of ${ A}_{\bar{0}}$ (resp. ${ A}_{\bar{1}}$, ${ A}_{\bar{2}}$) and is of grade 0 (resp. 1, 2). This gives a $Z_3$-graded structure to the algebra of matrices, in the sense that $\tau(T_iT_j)=\tau(T_i)+\tau(T_j)$ (mod 3). It is easy to check that the product of two $Z_3$-graded matrices is also a $Z_3$-graded matrix. As it can easily be shown, matrices of the form (\ref{3}) form a group provided that $ad-\beta\gamma\ne0$. We denote this group by $\widetilde{GL}(2)$.

\subsection{The algebra ${ O}(\tilde{M}_q(2))$}\label{sec3.1}

To determine a $q$-analogue of the algebra ${ O}(\tilde{M}(2))$, we will first  obtain the commutation relations between the matrix elements of the matrix $T$.

If ${ A}$ and ${ B}$ are $Z_3$-graded algebras, then their tensor product ${ A}\otimes{ B}$ is the $Z_3$-graded algebra whose underlying space is $Z_3$-graded tensor product of ${ A}$ and ${ B}$. The following definition \cite{Majid2} gives the product rule for tensor product of algebras.

\begin{definition}
If ${ A}$ is a $Z_3$-graded algebra, then the product rule in the $Z_3$-graded algebra ${ A}\otimes { A}$ is defined by
\begin{equation} \label{4}
(a_1\otimes a_2)(a_3\otimes a_4) = q^{\tau(a_2)\tau(a_3)} a_1a_3\otimes a_2a_4
\end{equation}
where $a_i$'s are homogeneous elements in the algebra ${ A}$.
\end{definition}

\noindent{\bf Remark 1.} It is well known that, the matrix $T$ given in (\ref{3}) defines the linear transformation
$T:\tilde{\mathbb R}_q^{0|2}\longrightarrow \tilde{\mathbb R}_q^{0|2}$ and $T:\tilde{\mathbb R}_q^{*0|2}\longrightarrow \tilde{\mathbb R}_q^{*0|2}$. As a result of these,  we have $T\Theta=\Theta' \in \tilde{\mathbb R}_q^{0|2}$ and $T\Phi=\Phi' \in \tilde{\mathbb R}_q^{*0|2}$, where $\Theta=(\theta, \varphi)^t$ and $\Phi=(\xi, x)^t$. However, the relation $\alpha_1\alpha_2 = q^{\tau(\alpha_1) \tau(\alpha_2)} \alpha_2\alpha_1$ for all elements $\alpha_1$ and $\alpha_2$ in the $Z_3$-graded algebra is inconsistent. Therefore, we will use the following transform while getting the commutation relations between the matrix elements of $T$.

Let $a$, $\beta$, $\gamma$, $d$ be elements of the algebra ${ O}(\tilde{M}(2))$. We also assume that the generators $a$ and $d$ are of degree 0, the generators $\gamma$ and $\beta$ are of degree 1 and 2, respectively. Then we can change the coordinates of a vector in $\tilde{\mathbb R}_q^{0|2}$ as follows
\begin{equation} \label{5}
\Theta' = \begin{pmatrix}\theta' \\ \varphi' \end{pmatrix} := \begin{pmatrix} a & \beta \\\gamma & d\end{pmatrix}\dot{\otimes}\begin{pmatrix}\theta\\ \varphi\end{pmatrix}, \quad
\Theta'' = \begin{pmatrix}\theta'' \\ \varphi'' \end{pmatrix} := \begin{pmatrix}\theta\quad \varphi\end{pmatrix}\dot{\otimes}\begin{pmatrix} a & \beta \\\gamma & d\end{pmatrix}.
\end{equation}
So, we can give the following proposition that can be proved with straightforward computations.

\begin{proposition}
The coordinates of $\Theta'$ and $\Theta''$ satisfy (\ref{1}) if and only if the generators $a,\beta,\gamma,d$ fulfill the relations
\begin{eqnarray} 
a \beta &=& \beta a, \quad \beta \gamma = \gamma \beta, \quad d\beta = \beta d, \label{6}\\
a \gamma &=& q \gamma a, \quad d\gamma = q^2 \gamma d, \label{7}\\
ad &=& da + (q-1) \beta \gamma, \label{8}
\end{eqnarray}
where $q$ is a cubic root of unity.
\end{proposition}

\noindent{\bf Remark 2.} Unlike the usual quantum group \cite{Manin1}, one interesting feature is that the element $\beta$ belongs to the center of the algebra.

\begin{definition} 
The $Z_3$-graded algebra ${ O}(\tilde{M}_q(2))$ is the quotient of the free algebra $k\{a,\beta,\gamma,d\}$ by the two-sided ideal $J_q$ generated by the six relations (\ref{6})-(\ref{8}) of Proposition 3.2.
\end{definition}

By relation (\ref{8}), we have
\begin{equation} \label{9}
{ D}_q := ad - q \beta \gamma = da - \beta \gamma.
\end{equation}
This element of ${ O}(\tilde{M}_q(2))$ is called the {\it $Z_3$-graded determinant}.

The proof of the following assertion is given by direct computation using the relations (\ref{6})-(\ref{8}).

\noindent{\bf Remark 3.} The $Z_3$-graded quantum determinant defined in (\ref{9}) commutes with $a$, $\beta$, $\gamma$ and $d$, so that the requirement ${ D}_q = 1$ is consistent.

\begin{proposition}
Let $T$ and $T'$ be two matrices such that their matrix elements satisfy the relations (\ref{6})-(\ref{8}). If all elements of $T$ commute according to the rule (\ref{4}) with all elements of $T'$, then the elements of the matrix (tensor) product $TT'$ obey the relations (\ref{6})-(\ref{8}). We also have
\begin{equation*}
{ D}_q(T\dot{\otimes}T') = { D}_q(T)\otimes { D}_q(T').
\end{equation*}
\end{proposition}

\noindent {\it Proof.} 
Let the matrix (tensor) product of $T$ with $T'$ be
\begin{equation*}
T\dot{\otimes}T' = \begin{pmatrix} a & \beta \\ \gamma & d\end{pmatrix}\dot{\otimes}\begin{pmatrix} a' & \beta' \\ \gamma' & d' \end{pmatrix} = \begin{pmatrix} X & Y \\ Z & W \end{pmatrix}.
\end{equation*}
Then using the relations (\ref{6})-(\ref{8}) with (\ref{4}) we get
\begin{eqnarray*}
XY
&=& (a\otimes a' + \beta\otimes \gamma')(a\otimes \beta' + \beta\otimes d') \\
&=& a^2\otimes a'\beta' + a\beta\otimes a'd' + \beta a\otimes \gamma'\beta' + q^2 \beta^2\otimes \gamma' d' \\
&=& a^2\otimes \beta' a' + a\beta \otimes d'a' + qa\beta\otimes \beta'\gamma' + \beta^2\otimes d'\gamma' \\
YX
&=& (a\otimes \beta' + \beta\otimes d')(a\otimes a' + \beta\otimes \gamma') \\
&=& a^2\otimes \beta' a' + a\beta \otimes d'a' + qa\beta\otimes \beta'\gamma' + \beta^2\otimes d'\gamma'.
\end{eqnarray*}
It can be similarly shown that relations $XZ=qZX$, $YZ=ZY$, etc., are provided. Proof of the latter as follows:
\begin{eqnarray*}
XW
&=& a\gamma\otimes a'\beta' + ad\otimes a'd' + q\beta \gamma\otimes \gamma'\beta' + \beta d\otimes \gamma' d' \\
YZ
&=& q^2 a\gamma\otimes \beta' a' + ad \otimes \beta'\gamma' + \beta\gamma\otimes d'a' + \beta d\otimes d'\gamma'\\
XW-qYZ
&=& ad\otimes (a'd' - q \beta'\gamma') - q \beta\gamma\otimes (d'a' - \gamma'\beta')
\end{eqnarray*}
and so ${ D}_q(T\dot{\otimes}T')$ reduces to ${ D}_q(T)\otimes { D}_q(T')$. \hfill$\square$

\subsection{Bialgebra structure on $\tilde{M}_q(2)$}\label{sec3.2}

We now supply the algebra ${ O}(\tilde{M}_q(2))$ with a bialgebra structure. The comultiplication and the counit will be the same as the usual quantum groups.

\begin{proposition}
$(1)$ There exist $Z_3$-graded algebra homomorphisms
\begin{equation*}
\Delta: { O}(\tilde{M}_q(2))\longrightarrow { O}(\tilde{M}_q(2))\otimes { O}(\tilde{M}_q(2)), \quad \epsilon: { O}(\tilde{M}_q(2))\longrightarrow {\mathbb C}
\end{equation*}
uniquely determined by
\begin{eqnarray} 
\Delta(a) = a\otimes a + \beta\otimes \gamma, \quad \Delta(\beta) = a\otimes \beta + \beta\otimes d, \label{10}\\
\Delta(\gamma) = \gamma\otimes a + d\otimes \gamma, \quad \Delta(d) = \gamma\otimes \beta + d\otimes d, \label{11}\\
\epsilon(a) = 1 = \epsilon(d), \quad \epsilon(\beta) = 0 = \epsilon(\gamma). \label{12}
\end{eqnarray}

\noindent $(2)$ With these maps, the algebra ${ O}(\tilde{M}_q(2))$ is a bialgebra which is neither commutative nor cocommutative.

\noindent $(3)$ The quantum determinant ${ D}_q$ is group-like element of ${ O}(\tilde{M}_q(2))$.
\end{proposition}

\noindent {\it Proof.} 
(1) In order to prove that $\Delta$ and $\epsilon$ are algebra homomorphisms, it is enough to show that the relations (\ref{6})-(\ref{8}) remain invariant under $\Delta$ and $\epsilon$. As an sample let us show that $\Delta(a\beta)=\Delta(\beta a)$:
\begin{eqnarray*}
\Delta(a\beta)
&=& \Delta(a)\Delta(\beta) = (a\otimes a + \beta\otimes \gamma)(a\otimes \beta + \beta\otimes d) \\
&=& a^2\otimes a\beta + a\beta\otimes ad + \beta a\otimes \gamma\beta + q^2 \beta^2\otimes \gamma d \\
&=& a^2\otimes \beta a + \beta a\otimes da + qa\beta\otimes \beta\gamma + \beta^2\otimes d\gamma \\
\Delta(\beta a)
&=& \Delta(\beta)\Delta(a) = (a\otimes \beta + \beta\otimes d)(a\otimes a + \beta\otimes \gamma) \\
&=& a^2\otimes \beta a + q a\beta \otimes \beta\gamma + \beta a\otimes da + \beta^2\otimes d\gamma.
\end{eqnarray*}
Analogously, one can prove another relations. For $\epsilon$ it is completely analogous.

\noindent (2) It is not difficult to check that the comultiplication $\Delta$ is coassociative in the sense that
\begin{equation} \label{13}
(\Delta \otimes \textrm{id}) \circ \Delta = (\textrm{id} \otimes \Delta) \circ \Delta
\end{equation}
and the counit $\epsilon$ has the property
\begin{equation} \label{14}
m \circ (\epsilon \otimes \mbox{id}) \circ \Delta = \mbox{id} = m \circ (\mbox{id} \otimes \epsilon) \circ \Delta.
\end{equation}
It follows that ${ O}(\tilde{M}_q(2))$ is indeed a bialgebra.

\noindent (3) To prove that the  $Z_3$-graded determinant ${ D}_q$ is group-like, it is enough to show that
\begin{equation} \label{15}
\Delta({ D}_q)={ D}_q \otimes { D}_q \quad \mbox{and} \quad \epsilon({ D}_q)=1.
\end{equation}
Indeed, some computations give
\begin{eqnarray*}
\Delta({ D}_q) &=& \Delta(a)\Delta(d)-q\Delta(\beta)\Delta(\gamma)\\
&=& ad\otimes ad + q\beta\gamma\otimes \beta\gamma - qad\otimes \beta\gamma - q\beta\gamma\otimes da \\
&=& ad\otimes (ad-q\beta\gamma) + q \beta\gamma\otimes (\beta\gamma-da) \\
&=& (ad-q\beta\gamma)\otimes (da-\beta\gamma)
\end{eqnarray*}
and $\epsilon(ad-q\beta\gamma)=\epsilon(a)\epsilon(d)-q\epsilon(\beta)\epsilon(\gamma)=1$. \hfill$\square$

The bialgebra ${ O}(\tilde{M}_q(2))$ is called the coordinate algebra of the $Z_3$-graded (quantum) matrix space $\tilde{M}_q(2)$.

\subsection{The $Z_3$-graded Hopf algebra ${ O}(\widetilde{GL}_q(2))$}\label{sec3.3}

Using the quantum determinant ${ D}_q$ belonging to  the algebra ${ O}(\tilde{M}_q(2))$, we can define a {\it new} Hopf algebra adding an inverse $t^{-1}$ to ${ O}(\tilde{M}_q(2))$. Let ${ O}(\widetilde{GL}_q(2))$ be the quotient of the algebra ${ O}(\tilde{M}_q(2))$ by the two-sided ideal generated by the element $t{ D}_q-1$. For short we write
\begin{equation*}
{ O}(\widetilde{GL}_q(2)) := { O}(\tilde{M}_q(2))[t]/\langle t{ D}_q-1\rangle.
\end{equation*}
Then the algebra ${ O}(\widetilde{GL}_q(2))$ is again a bialgebra.

\begin{lemma}
The elements of the matrix
\begin{equation} \label{16}
\tilde{T} = \begin{pmatrix} \tilde{a} & \tilde{\beta} \\ \tilde{\gamma} & \tilde{d} \end{pmatrix} = \begin{pmatrix} d & -\beta \\ -q\gamma & a \end{pmatrix}
\end{equation}
satisfy the defining relations of the algebra ${ O}(\widetilde{GL}_{q^2}(2))$ and thus ${ O}(\widetilde{GL}_{q^2}(2))$ is the opposite algebra of
${ O}(\widetilde{GL}_q(2))$.
\end{lemma}

\noindent {\it Proof.} 
The use of relations (\ref{6})-(\ref{8}) imply
\begin{eqnarray*}
\tilde{a} \tilde{\beta} &=& \tilde{\beta} \tilde{a}, \quad \tilde{\beta} \tilde{\gamma} = \tilde{\gamma} \tilde{\beta}, \quad \tilde{\beta} \tilde{d} = \tilde{d} \tilde{\beta}, \\
\tilde{a} \tilde{\gamma} &=& q^2 \tilde{\gamma} \tilde{a}, \quad \tilde{\gamma} \tilde{d} = q^2 \tilde{d} \tilde{\gamma}, \\
\tilde{a} \tilde{d} &=& \tilde{d} \tilde{a} + (1-q^2) \tilde{\beta} \tilde{\gamma},
\end{eqnarray*}
which are the defining relations of the algebra ${ O}(\widetilde{GL}_{q^2}(2))$. The second claim follows from the fact that $q^3=1$. \hfill$\square$
 
\begin{proposition}
The bialgebra ${ O}(\widetilde{GL}_q(2))$ is a $Z_3$-graded Hopf algebra. The antipode $S$ of ${ O}(\widetilde{GL}_q(2))$ is given by
\begin{equation} \label{17}
S(a) = d \,{ D}_q^{-1}, \quad S(\beta) = - \beta \,{ D}_q^{-1}, \quad S(\gamma) = -q\gamma \,{ D}_q^{-1}, \quad S(d) = a \,{ D}_q^{-1}.
\end{equation}
\end{proposition}

\noindent {\it Proof.} 
By Lemma 3.6, there exists an algebra anti-homomorphism $S$ from ${ O}(\widetilde{GL}_q(2))$ to ${ O}(\widetilde{GL}_{q^2}(2))$ such that $S(a)=\tilde{a}$, etc.
To prove that $S$ is an antipode for ${ O}(\widetilde{GL}_q(2))$, we have to check the antipode axiom
\begin{equation} \label{18}
m \circ (S \otimes \mbox{id}) \circ \Delta = \epsilon = m \circ (\mbox{id} \otimes S) \circ \Delta
\end{equation}
for the generators. To check the axiom (\ref{18}) for the generators is equivalent to verify the following matrix equality
\begin{equation*}
T \tilde{T} { D}_q = \epsilon(T) = \tilde{T} T { D}_q
\end{equation*}
which follows from ${ D}_q = ad - q \beta \gamma$ in ${ O}(\widetilde{GL}_q(2))$ with $S(T)={ D}_q^{-1} \tilde{T}=T^{-1}$. The details can be checked easily. \hfill$\square$

\begin{definition} 
The $Z_3$-graded Hopf algebra ${ O}(\widetilde{GL}_q(2))$ is called the coordinate algebra of the $Z_3$-graded (quantum) group $\widetilde{GL}_q(2)$.
\end{definition}

\subsection{Coactions on the $Z_3$-graded exterior plane}\label{sec3.4}

In bialgebra terminology, the second suggestion of Proposition 3.2 yields the following.
\begin{proposition}
The algebra ${ O}(\tilde{\mathbb R}_q^{0|2})$ is a left and right comodule algebra of the bialgebra ${ O}(\tilde{M}_q(2))$ with left coaction
$\delta_L$ and right coaction $\delta_R$ such that
\begin{eqnarray} 
\delta_L(\theta) = a \otimes \theta + \beta \otimes \varphi, \quad \delta_L(\varphi) = \gamma \otimes \theta + d \otimes \varphi, \label{19}\\
\delta_R(\theta) = \theta \otimes a + \varphi \otimes \gamma, \quad \delta_R(\varphi) = \theta \otimes \beta + \varphi \otimes d. \label{20}
\end{eqnarray}
\end{proposition}

\noindent {\it Proof.} 
It is not difficult to verify that (\ref{19}) and (\ref{20}) define algebra homomorphisms $\delta_L$ from $O(\tilde{\mathbb R}_q^{0|2})$ to
${ O}(\tilde{M}_q(2))\otimes { O}(\tilde{\mathbb R}_q^{0|2})$ and $\delta_R$ from ${ O}(\tilde{\mathbb R}_q^{0|2})$ to
${ O}(\tilde{\mathbb R}_q^{0|2})\otimes { O}(\tilde{M}_q(2))$, respectively. It remains to be checked that $\delta_L$ and $\delta_R$ are coactions, i.e., the conditions
\begin{equation} \label{21}
(\Delta\otimes\textrm{id})\circ\delta_L = (\textrm{id}\otimes\delta_L)\circ\delta_L, \quad m\circ(\epsilon\otimes\mbox{id})\circ\delta_L = \mbox{id} \end{equation}
and
\begin{equation} \label{22}
(\textrm{id}\otimes\Delta)\circ\delta_R = (\delta_R\otimes\textrm{id})\circ\delta_R, \quad m\circ(\mbox{id}\otimes\epsilon)\circ\delta_R = \mbox{id} \end{equation}
are satisfied. For examples,
\begin{eqnarray*}
(\Delta\otimes\textrm{id})\delta_L(\theta)
&=& (\Delta\otimes\textrm{id})(a \otimes \theta + \beta \otimes \varphi)\\
&=& (a\otimes a + \beta\otimes \gamma)\otimes \theta + (a\otimes \beta + \beta\otimes d)\otimes \varphi \\
&=& a \otimes(a\otimes \theta + \beta \otimes \varphi) + \beta \otimes(\gamma \otimes \theta + d \otimes \varphi) \\
&=& a \otimes \delta_L(\theta) + \beta \otimes \delta_L(\varphi)\\
&=& (\textrm{id}\otimes\delta_L)\delta_L(\theta)
\end{eqnarray*}
and
\begin{eqnarray*}
m\circ(\epsilon\otimes\mbox{id})\delta_L(\theta)
&=& m(\epsilon\otimes\mbox{id})(a \otimes \theta + \beta \otimes \varphi)\\
&=& m(1\otimes \theta + 0\otimes \varphi) \\
&=& \theta
\end{eqnarray*}
as expected. \hfill$\square$

\noindent{\bf Remark 4.} In fact, there exists a left coaction of ${ O}(\tilde{\mathbb R}_q^{*0|2})$ on the plane $\tilde{\mathbb R}_q^{*0|2}$, called a left comodule-${ O}(\tilde{\mathbb R}_q^{*0|2})$ satisfying the conditions (\ref{21}).

\noindent{\bf Remark 5.} An easy computation shows that the ideal $(\vartheta:=\theta\varphi-q^2\theta\varphi)$ of $\tilde{\mathbb R}_q^{0|2}$ is a subcomodule
of $\tilde{\mathbb R}_q^{0|2}$. The proof is immediate: Indeed, since $\delta_L$ is an algebra map, it is only necessary to show that $\delta_L(\vartheta)={ D}_q\otimes \vartheta$. Using relations (\ref{6})-(\ref{8}) with (\ref{9}) we get
\begin{eqnarray*}
\delta_L(\vartheta) &=& \delta_L(\theta) \delta_L(\varphi) - q^2 \delta_L(\varphi) \delta_L(\theta)\\
&=& qa\gamma \otimes \theta^2 + ad \otimes \theta\varphi + q^2 \beta\gamma \otimes \varphi\theta + \beta d \otimes \varphi^2 - q^2 \gamma a \otimes \theta^2 \\
&& - q \gamma\beta \otimes \theta\varphi - q^2 da \otimes \varphi\theta - d\beta \otimes \varphi^2\\
&=& (ad - q \beta \gamma) \otimes \theta\varphi - q^2 (da - \beta \gamma) \otimes \varphi\theta = { D}_q \otimes \vartheta
\end{eqnarray*}
as expected. \hfill$\square$

\subsection{The Hopf algebra ${ O}(\widetilde{SL}_q(2))$}\label{sec3.5}

We know that, since the determinant ${ D}_q$ is group-like, the two-sided ideal $\langle { D}_q-1\rangle$ generated by the element ${ D}_q-1$ is a biideal of ${ O}(\widetilde{M}_q(2))$. So the quotient ${ O}(\widetilde{SL}_q(2)) := { O}(\tilde{M}_q(2))/\langle{ D}_q-1\rangle$ is a bialgebra.

\begin{proposition}
There exists a Hopf $\star$-algebra structures on the Hopf algebra ${O}(\widetilde{SL}_q(2))$ such that
\begin{equation} \label{23}
a^\star = a, \quad \beta^\star = \beta, \quad \gamma^\star = q \gamma, \quad d^\star =d.
\end{equation}
\end{proposition}


\section{$Z_3$-graded quantum algebra of $\widetilde{GL}_q(2)$}\label{sec4}

In this section, using the method of \cite{FRT}, we give an $R$-matrix formulation for the $Z_3$-graded quantum group $\widetilde{GL}_q(2)$ and obtain a $Z_3$-graded universal enveloping algebra $U_q(\widetilde{gl}(2))$.

\subsection{The FRT construction for $\widetilde{GL}_q(2)$}\label{sec4.1}

The $R$-matrix formulation (the FRT-relation $\hat{R} T_1T_2 = T_1T_2 \hat{R}$) for the quantum matrix groups \cite{FRT} can be considered as a compact matrix form of the commutation relations between the generators of an associative algebra.

The formulation for the $Z_3$-graded quantum group $\widetilde{GL}_q(2)$ has the same form, but matrix tensor product includes additional $q$-factors related
to $Z_3$-grading. Two matrices $A$, $B$ $(\tau(A_{ij}) = \tau(i) + \tau(j))$ are multiplied according to the rule
\begin{equation} \label{24}
(A\otimes B)_{ij,kl} = q^{\tau(j)(\tau(i)+\tau(k))} A_{ik}B_{jl}.
\end{equation}
Due to this prescription, $T_2 = I\otimes T$ has the same block-diagonal form as in the standard (ungraded) case while $T_1 = T\otimes I$ includes the additional factors $q$ for graded elements standing at some of odd rows of blocks. For the $Z_3$-graded quantum group $\widetilde{GL}_q(2)$ the $R$-matrix satisfying the $Z_3$-graded Yang-Baxter equation has in the form
\begin{eqnarray} \label{25}
\hat{R} = \begin{pmatrix}
q & 0 & 0 & 0 \\
0 & 0 & 1 & 0 \\
0 & 1 & q-q^2 & 0 \\
0 & 0 & 0 & q
\end{pmatrix}=(\hat{R}_{ij}^{kl})
\end{eqnarray}
where $\hat{R}=\underline{P}R$ and $\underline{P}$ denotes the $Z_3$-graded permutation operator defined by
$\underline{P}(a\otimes b)=q^{\tau(a)\tau(b)} b\otimes a$ on homogeneous elements. A simple calculation shows that this operator represents the 3rd-root of the permutation operator $P$ with action $P(a\otimes b)=b\otimes a$.

The condition for the matrices to belong to the $Z_3$-graded quantum group $\widetilde{GL}_q(2)$ is given below, but it will not be proved here.

\begin{proposition}
A 2x2-matrix $T$ is a $Z_3$-graded quantum matrix if and only if
\begin{equation} \label{26}
\hat{R} T_1T_2 = T_1T_2 \hat{R}
\end{equation}
where matrix elements of $T$ are $Z_3$-graded.
\end{proposition}

\subsection{A $Z_3$-graded universal enveloping algebra $U_q(\widetilde{gl}(2))$}\label{sec4.2}

The $Z_3$-graded quantum algebra of $\widetilde{GL}_q(2)$ can be analogous construction to approach of the Leningrad school. The $Z_3$-graded quantum algebra of $\widetilde{GL}_q(2)$ has four generators: $U$ and $V$ are of degree 0, $X_-$ and $X_+$ are of degrees 1 and 2, respectively.

\begin{proposition}
The generators of the $Z_3$-graded quantum algebra satisfy the following relations
\begin{eqnarray} 
&& UV=VU, \quad UX_\pm = q^{\pm2} X_\pm U, \quad VX_\pm = q^{\mp2} X_\pm V, \label{27}\\
&& X_+X_- - X_-X_+ = \frac{UV^{-1} - VU^{-1}}{q^2-q} \label{28}
\end{eqnarray}
\end{proposition}

\noindent {\it proof.} 
The generators $U$, $V$, $X_\pm$ can be written in two 2x2 matrix as follows
\begin{eqnarray} \label{29}
L^+ = \begin{pmatrix} U & \lambda X_+ \\ 0 & V \end{pmatrix}, \quad L^- = \begin{pmatrix} U^{-1} & 0 \\ \lambda X_- & V^{-1} \end{pmatrix}
\end{eqnarray}
where $\lambda=q-q^2$. The  matrices $L^\pm$ satisfy the following relations
\begin{equation} \label{30}
R^+ L^\pm_1 L^\pm_2 = L^\pm_2 L^\pm_1 R^+,
\end{equation}
where the matrix $R^+$ is defined by $R^+=\underline{P}R\underline{P}$. The relations (\ref{27}) follow from the relations (\ref{30}).
To obtain the relation (\ref{28}) we use the relation
\begin{equation} \label{31}
R^+ L^-_1 L^+_2 = L^+_2 L^-_1 R^+. 
\end{equation}

\begin{proposition}
The coproduct of the generators is given by
\begin{equation} \label{32}
\Delta(L^\pm) = L^\pm \dot{\otimes} L^\pm.
\end{equation}
\end{proposition}



\baselineskip=10pt
\begin{thebibliography}{00}
\bibitem{Manin1} Yu I. Manin, {\it Quantum groups and noncommutative geometry}, Montreal Univ. Preprint, 1988.
\bibitem{Manin2} Yu I. Manin, {\it Multiparametric quantum deformation of the general linear supergroup}, Commun. Math. Phys. 123 (1989) 163-175.
\bibitem{Chung} W. S. Chung, {\it Quantum Z$_3$-graded space}, J. Math. Phys. {\bf 35} (1994) 2497-2504.
\bibitem{Celik} S. Celik, {\it Differential geometry of the Z$_3$-graded quantum superplane}, J. Phys. A: Math. Gen. {\bf 35} (2001) 4257-4268.
\bibitem{Kli-sch} A. Klimyk and K. Schm\"udgen, {\it Quantum Groups and Their Representations}, Texts and Monographs in Physics, Springer, New York et al., 1997.
\bibitem{Majid1} Majid, S, {\it Foundations of Quantum Group Theory}, Cambridge Univ. Press, Cambridge, 1995.
\bibitem{Majid2} S. Majid, {\it Anyonic quantum groups}, In Spinors, Twistors, Clifford Algebras and Quantum Deformations (Proc. of 2nd Max Born Symposium, Wroclaw, Poland, 1992), Z. Oziewicz et al, eds., pages 327-336.
\bibitem{FRT} L. D. Faddeev, N. Yu. Reshetikhin and L. A. Takhtajan, {\it Quantization of Lie groups and Lie algebras}, Leningrad Math. J. {\bf 1} (1990) 193-225.
\end{thebibliography}
\end{document}